\documentclass[11pt]{article}
\usepackage{amsmath, amssymb, mathrsfs}
\usepackage{graphicx}
\usepackage{geometry}
\usepackage{graphicx}
\usepackage{wrapfig}
\usepackage{caption}

\begin{document}

\begin{center}

\Large
\textbf{On the occasion of Dr. Ivan Dmitrievich Remizov's 40th birthday}

\large

\vskip5mm

Ksenia Dragunova\footnote{HSE University}, 
Elmira Kalimulina\footnote{Institute for Information Transmission Problems of the Russian Academy of Sciences}, 
Rustem Kalmetev\footnote{Moscow Institute of Physics and Technology}, 
Nasrin Nikbakht\footnote{University of Auckland,  corresponding author, email nasrin.nikbakht@gmail.com}, 
Elena Nozdrinova\footnote{HSE University}, 
Alexander Pechen\footnote{Steklov Mathematical Institute of Russian Academy of Sciences}, 
Alexei Savvateev\footnote{Moscow Institute of Physics and Technology}, 
Evgeni Shavgulidze\footnote{Lomonosov Moscow State University}, 
Ranjit Singh\footnote{Peoples' Friendship University of Russia},
Alexander Vedenin\footnote{HSE University}, 
Boris Volkov\footnote{Moscow Institute of Physics and Technology}

\end{center}

\normalsize

\begin{abstract}
This article celebrates the 40th anniversary of Dr. Ivan Dmitrievich Remizov, a mathematician who made a number of important contributions to the theory of one-parameter operator semigroups -- a branch of functional analysis which has applications to differential equations, mathematical physics, random processes, control theory, and quantum mechanics. Born on December 7, 1984, Dr. Remizov obtained his Ph.D. in 2018 from Moscow State University and has made substantial contributions, particularly in the study of Chernoff approximations of one-parameter semigroups of operators. This article reviews his academic background, research achievements, and his impact on the mathematical community.
\end{abstract}

\section{Life and career overview}

\begin{wrapfigure}{l}{0.25\textwidth}
    \centering    
\includegraphics[width=0.25\textwidth]{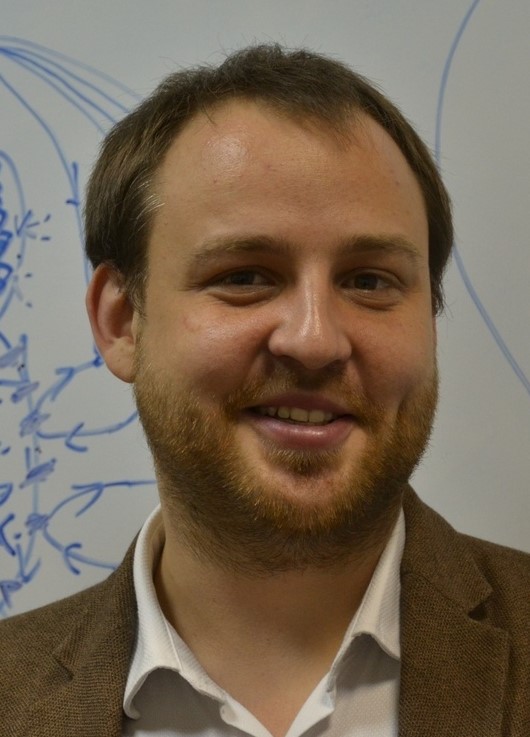}
\caption*{Photo by E.~Nozdrinova, N.~Novgorod city, 2019}
\end{wrapfigure}

Dr. Ivan Dmitrievich Remizov approached his 40th birthday on December 7, 2024, and it is our pleasure to commemorate his influential career. Born on December 7, 1984, in Gorky, USSR (now Nizhny Novgorod, Russia), Ivan Remizov atended meetings of the Nizhny Novgorod Mathematical Society since 1999, and was officially accepted as a regular member of the Society in 2001. He graduated from physics and mathematics focused Lyceum No 40 in 2002 and became a student of the N.Novgorod State University, department of Mathematics and Mechanics, where he received first lessons in Functional Analysis from Dr. Alexander Abrosimov \cite{Abrosimov}. Willing to continue studies in that field, Remizov in 2005 transferred to the Moscow State University (Department of Mathematics and Mechanics), where he became a student of Professor Oleg Georgievich Smolyanov \cite{Smolyanov}. At MSU Remizov obtained graduate diploma in 2008 and Ph.D. diploma in 2018 under the supervision of Prof. O.G.Smolyanov. In 2010-2014 Remizov was a teaching assistant in the MSU, and in 2014-2017 he was an Assistant Professor at the Bauman Moscow State Technical University. In 2018 Remizov was invited by Prof. O.V.Pochinka \cite{Pochinka} to take the position of Associate Professor at the N.Novgorod branch of the HSE University, so Remizov moved back to Nizhny Novgorod and was teaching at HSE NN until 2023. In 2023 Remizov moves back to Moscow and takes a postdoctoral position under supervision of Prof. A.A.Shkalikov \cite{Shkalikov} at the MSU, also keeping the position of a Senior Research Fellow at the HSE NN (distant work). Since 2025 Remizov has also held a Senior Research Fellow position at the IITP Dobrushin Math. Laboratory (part-time work).

The mathematical genealogy of Dr. Remizov is as follows. His Ph.D. supervisor was O.G.Smolyanov, who defended a Ph.D. thesis under the supervision of G.E. Shilov, who defended a Ph.D. thesis under the supervision of I.M. Gelfand, who defended a Ph.D. thesis under the supervision of A.N. Kolmogorov. It is worth mentioning that O.G.Smolyanov considered S.V.Fomin as his main teacher, and S.V.Fomin defended his Ph.D. thesis under the supervision of A.N. Kolmogorov. During his career path, Ivan Remizov has established himself as a strong mathematician in functional analysis and its applications, a bright representative of the Moscow School of Function Theory.

\section{Research Field}

Dr. Remizov’s research is focused on the theory of strongly continuous one-parameter semigroups of linear bounded operators in a Banach space, also called $C_0$-semigroups. $C_0$-semigroup is a generalization of the notion of an exponent \cite{EN1, Gold}. If $\mathcal{F}$ is a Banach space, and $L$ is a densely defined closed linear operator in $\mathcal{F}$, then $V$ is called a $C_0$-semigroup with generator $L$ iff the mapping $V\colon [0,+\infty)\to\mathscr{L}(\mathcal{F})$ is given (i.e. for each  $t\geq0$ operator $V(t)$ maps $\mathcal{F}$ to $\mathcal{F}$ linearly and continuously), $V(0)=I$, $V(t_1+t_2)=V(t_1)V(t_2)$ for all $t_1,t_2\geq 0$, function $t\mapsto V(t)f$ is continuous for each $f\in\mathcal{F}$, $V'(0)=L$. The notation $V(t)=e^{tL}$ is used. The definition of $C_0$-semigroup does not provide a direct method of finding $e^{tL}$ if one knows $L$. If $L$ is a bounded operator, a finite matrix, or a number, then the standard power series $\sum_{k=0}^\infty(tL)^k/k!=e^{tL}$ can be used for this purpose, but if $L$ is an unbounded operator, then this series diverges. 

In 1968 Paul Chernoff published \cite{Chernoff} his theorem, which allows to find $e^{tL}$ approximately if a so-called Chernoff function is known. One of the formulations of the Chernoff theorem says that if 1) the operator-valued function $C\colon [0,+\infty)\to\mathscr{L}(\mathcal{F})$ is given, $C(0)=I$, the function $t\mapsto C(t)f$ is continuous for each $f\in\mathcal{F}$, $C'(0)\subset L$, $\overline{C'(0)}=L$ (these assumptions are called conditions of Chernoff tangency of $C$ to $L$ according to Remizov's terminology) and 2) there exists $w\geq 0$ such that $\|C(t)\|\leq e^{tw}$ for each $t\geq 0$, then $e^{tL}f=\lim_{n\to\infty} C(t/n)^nf$ for each $f\in\mathcal{F}$. In that case $C$ is called a Chernoff function for operator $L$, and $C(t/n)^n$ are called Chernoff approximations for $e^{tL}$. In the one-dimensional case $\mathcal{F}=\mathscr{L}(\mathcal{F})=\mathbb{R}$ the Chernoff theorem says that if $c\in C([0,+\infty),\mathbb{R})$, $c(0)=1$, $c'(0)=l$ then $c(t/n)^n=(1+ tl/n +o(1/n))^n\to e^{tl}$ which is a simple fact from calculus. 

In the general case where $L$ is an unbounded operator in infinite-dimensional Banach space $\mathcal{F}$ the Chernoff approximation method is often the only practical method of finding the exponent $e^{tL}$. In particular, if $L$ is a differential operator with variable coefficients then (in many cases, it was shown by O.G.Smolyanov and his students and colleagues, including I.D.Remizov) it is possible to find a Chernoff function $C$ for $L$ and express it by a short formula which explicitly includes the variable coefficients of $L$. Thus, Chernoff approximations $C(t/n)^n$ that are arbitrary close to $e^{tL}$ are expressed explicitly in terms of variable coefficients of operator $L$, where the coefficients run over some class of functions and play the role of parameters of the problem. There are also methods for finding $e^{tL}$ based on the resolvent of $L$, but the resolvent $(\lambda -L)^{-1}$, as well as the exponent $e^{tL}$, is defined by $L$ uniquely and thus is difficult to find, while the Chernoff function is not defined by $L$ uniquely, so there is some freedom in finding Chernoff functions for a given $L$. This makes the method of Chernoff approximation a flexible and powerful tool of functional analysis with various applications \cite{Butko-2019}.

If the Chernoff function $C$ for operator $L$ is known, then the solution of the Cauchy problem $[U'(t)=LU(t),U(0)=u_0]$ for function $U\in C^1([0,+\infty),\mathcal{F})$ can be found as $U(t)=e^{tL}u_0=\lim_{n\to\infty}C(t/n)^nu_0$. If $\mathcal{F}$ is a Banach space of number-valued functions of variable $x$, then one can introduce the notation $u(t,x)=(U(t))(x)$ and rewrite this Cauchy problem as $[u_t(t,x)=Lu(t,x),u(0,x)=u_0(x)]$, and the solution of it is still given by the formula $u(t,x)=(e^{tL}u_0)(x)=\lim_{n\to\infty}(C(t/n)^nu_0)(x)$. If $L$ is the operator of taking the first derivative $(Lf)(x)=f'(x)$, then $u_t(t,x)=Lu(t,x)$ is the transport equation $u_t(t,x)=u_x(t,x)$. If $(Lf)(x)=f''(x)$, then $u_t(t,x)=Lu(t,x)$ is the heat equation $u_t(t,x)=u_{xx}(t,x)$ etc. This connection between $C_0$-semigroups and evolution equations is one of the reasons why $C_0$-semigroups and their Chernoff approximations have various applications to differential equations, mathematical physics, random processes, control theory and quantum mechanics \cite{EN1, Gold, Butko-2019}. If the Chernoff function $C$ is based on the integral operator, then $C(t/n)^n$ is an $n$-tuple integral with multiplicity tending to infinity, and its limit as $n\to\infty$ is the Feynman path integral. This is why in that case the representation $u(t,x)=\lim_{n\to\infty}(C(t/n)^nu_0)(x)$ is called a Feynman formula and can be used to calculate Feynman path integrals \cite{Maz, BGS2010, SS, BogS}.

\section{Research Contribution}

For several concrete classes of second-order differential operators $L$ Dr. Remizov has found Chernoff functions which were represented as elementary functions of variable coefficients of the operator $L$, and obtained Chernoff approximations to the solution of the Cauchy problem for the parabolic equation $[u_t(t,x)=Lu(t,x),u(0,x)=u_0(x)]$. For the case where $x$ runs over a separable Hilbert space, he found a Chernoff function based on the integral operator \cite{R1}. For the case where $x$ runs over $\mathbb{R}^1$ he proposed a new class of Chernoff functions based on the shift operator \cite{R2}, and then applied the same technique to the case $x\in\mathbb{R}^d$ for arbitrary $d\in\mathbb{N}$ \cite{R3}. Later, this approach was generalized by S.Mazzucchi and V.Moretti to the case of manifolds of bounded geometry and published in a joint paper with I.D.Remizov and O.G.Smolyanov \cite{R4}.

Dr. Remizov also proposed \cite{R5} a universal method of constructing Chernoff approximations for solutions of the Schodinger-type equations with Hamiltonian $H$, i.e. equations $u_t(t,x)=Lu(t,x)$ with $L=-iH$ where $H$ is a self-adjoint operator in Hilbert space. I.D.Remizov proved that if $S$ is Chernoff tangent to $H$, $S(t)=S(t)^*$, $0\neq a\in\mathbb{R}$, then one of Chernoff functions for the operator $iatH$ is given by Remizov's formula formula $R(t) = \exp(ia(S(t) - I))$. Here in the exponent only bounded operators appear, so this exponent can always be defined by a power series, meanwhile if $H$ is an unbounded operator, then the power series for $\exp(iatH)$ diverges. This allows finding $\exp(-itH)$ if one knows $S(t)=\exp(-tH)$ or $S(t)=\exp(tH)$ or even less -- a function $S$ that is Chernoff tangent to $H$ and satisfies the condition $S(t)=S(t)^*$. If $S$ is based on an integral operator, then the solution of the equation is represented as a sum of multiple integrals with multiplicity tending to infinity, such expressions were at first obtained by I.D.Remizov and were called \textit{quasi-Feynman formulas} by him. Using this technique and employing shift-based Chernoff functions, he found Chernoff approximations to the solution of the Cauchy problem for a standard Schr\"odinger equation in $\mathbb{R}^d$ and then to Schr\"odinger equation with variable coefficients and derivatives of arbitrary high order in $\mathbb{R}^1$ \cite{R6}.

Dr. Remizov studied the rate of convergence of Chernoff approximations and introduced \cite{R7} the notion of the \textit{approximation subspace} in the Chernoff theorem, i.e. the set of all vectors $f\in\mathcal{F}$ such that $\|e^{tL}f-C(t/n)^nf\|=o(a_n)$ for a given sequence $a_n\to0$. Later he proved that the rate of convergence in general can be arbitrarily high and arbitrarily slow, and together with O.E.Galkin he found conditions \cite{R8} on $C$ and $f$ which guarantee that $\|e^{tL}f-C(t/n)^nf\|=o(1/n^k)$ for given $k\in\mathbb{N}$. Together with his students K.A.Dragunova and N.Nikbakht for cases not covered by conditions of the Galkin-Remizov theorem he studied the rate of convergence numerically \cite{R9}. For approximations that converge faster than $1/n$ he proposed the term \textit{fast convergent} and for approximations that converge faster than $1/n^k$ for all $k\in\mathbb{N}$ he proposed the term \textit{super fast convergent}. Together with his student A.V.Vedenin for operator $Lf = af'' + bf' + cf$ where functions $a,b,c$ are uniformly continuous and bounded Dr. Remizov proposed a method of constructing fast converging Chernoff approximations, see \cite{V} and joint paper to appear. Also, Dr. Remizov posed an open problem: does practically valuable example of super fast converging Chernoff approximations to the operator $Lf = af'' + bf' + cf$ exist or not. In the case of positive answer, one can expect further increase of interest to the method of Chernoff approximation.

Using the fact that the resolvent of the generator of the $C_0$-semigoup is the Laplace transform of the semigroup $(\lambda -L)^{-1}=\int_0^{+\infty}e^{-\lambda t}e^{tL}dt$, Dr. Remizov proposed a method of construction of Chernoff approximations of solutions of a wide class of second order linear ODEs and elliptic PDEs with variable coefficients \cite{R10}. The idea of this approach is based on the fact that these equations can be written in the form $(\lambda -L)f=g$ where the variable coefficients of the operator $L$, as well as the function $g$, play the role of parameters, and $f$ is the function that we aim to find.

\section{Contribution to the Mathematical Community}

Dr. Remizov's influence extends to his active role in organizing and supporting mathematical initiatives. He started in 2021 (and chaired the Organizing Committee in 2021-2023) the series of online international conferences on one-parameter semigroups of operators (OPSO) and co-edited a list of open problems in this area. He initiated the “Review of Diploma Works of Russian Mathematicians” competition, serving as chair from 2019 to 2023, which encouraged academic excellence and collaboration among emerging mathematicians. Also, he helped to open an international Master program in Mathematics at HSE Nizhny Novgorod in 2018 and was teaching there during 2018-2023, also served as examiner for international students in the admission process.

\section{Conclusion}

As Dr. Ivan Remizov celebrates his 40th anniversary, we recognize his important contributions in analysis. His innovative research, leadership, and dedication to the mathematical community continue to inspire scholars and students. This milestone marks not only the accomplishments of a distinguished career, but also the promising future contributions Dr. Remizov is sure to make in advancing mathematical knowledge and collaboration. It is worth to mention that Dr. Remizov's interests are not only in the field of mathematics. Since 2022, he has been conducting research in psychophysiology and also practices as a psychotherapist. We wish Dr. Remizov good luck on his path.

\end{document}